# BOUNDS FOR BAYESIAN ORDER IDENTIFICATION WITH APPLICATION TO MIXTURES

By Antoine Chambaz and Judith Rousseau

*Université Paris Descartes and Université Dauphine*

The efficiency of two Bayesian order estimators is studied. By using nonparametric techniques, we prove new underestimation and overestimation bounds. The results apply to various models, including mixture models. In this case, the errors are shown to be $O(e^{-an})$ and $O((\log n)^b/\sqrt{n})$ $(a, b > 0)$, respectively.

**1. Introduction.** Order identification deals with the estimation and test of a structural parameter which indexes the complexity of a model. In other words, the most economical representation of a random phenomenon is sought. This problem is encountered in many situations, including: mixture models [13, 19] with an unknown number of components; cluster analysis [9], when the number of clusters is unknown; autoregressive models [1], when the process memory is not known.

This paper is devoted to the study of two Bayesian estimators of the order of a model. Frequentist properties of efficiency are particularly investigated. We obtain new efficiency bounds under mild assumptions, providing a theoretical answer to the questions raised, for instance, in [7] (see their Section 4).

1.1. *Description of the problem.* We observe $n$ i.i.d. random variables (r.v.) $(Z_1, \dots, Z_n) = Z^n$ with values in a measured sample space $(\mathcal{Z}, \mathcal{F}, \mu)$.

Let $(\Theta_k)_{k\geq 1}$ be an increasing family of nested parametric sets and $d$ the Euclidean distance on each. The dimension of $\Theta_k$ is denoted by $D(k)$. Let $\Theta_\infty = \bigcup_{k\geq 1} \Theta_k$ and for every $\theta \in \Theta_\infty$, let $f_\theta$ be the density of the probability measure $P_\theta$ with respect to the measure $\mu$.

The order of any distribution $P_{\theta_0}$ is the unique integer $k$ such that $P_{\theta_0} \in \{P_\theta : \theta \in \Theta_k \setminus \Theta_{k-1}\}$ (with convention $\Theta_0 = \varnothing$). It is assumed that the distribution $P^\star$ of $Z_1$ belongs to $\{P_\theta : \theta \in \Theta_\infty\}$. The density of $P^\star$ is denoted by $f^\star = f_{\theta^\star}$ ($\theta^\star \in \Theta_{k^\star} \setminus \Theta_{k^\star-1}$). The order of $P^\star$ is denoted by $k^\star$, and is the quantity of interest here.

We are interested in frequentist properties of two Bayesian estimates of $k^\star$. In that perspective, the problem can be restated as an issue of composite hypotheses testing (see [4]), where the quantities of interest are $P^\star\{\widetilde{k}_n < k^\star\}$ and $P^\star\{\widetilde{k}_n > k^\star\}$, the under- and over-estimation errors, respectively. In this paper we determine upper-bounds on both errors on $\widetilde{k}_n$ defined as follows.

Let $\Pi$ be a prior on $\Theta_\infty$ that writes as $d\Pi(\theta) = \pi(k)\pi_k(\theta)\,d\theta$, for all $\theta \in \Theta_k$ and $k \geq 1$. We denote by $\Pi(k|Z^n)$ the posterior probability of each $k \geq 1$. In a Bayesian decision theoretic perspective, the Bayes estimator associated with the 0–1 loss function is the mode of the posterior distribution of the order $k$:

$$\widehat{k}_n^{\mathrm{G}} = \arg\max_{k \geq 1}\{\Pi(k|Z^n)\}.$$

It is a global estimator. Following a more local and sequential approach, we propose another estimator:

$$\widehat{k}_n^{\mathrm{L}} = \inf\{k \geq 1 : \Pi(k|Z^n) \geq \Pi(k+1|Z^n)\} \leq \widehat{k}_n^{\mathrm{G}}.$$

If the posterior distribution on $k$ is unimodal, then obviously both estimators are equal. The advantage of $\widehat{k}_n^{\mathrm{L}}$ over $\widehat{k}_n^{\mathrm{G}}$ is that $\widehat{k}_n^{\mathrm{L}}$ does not require the computation of the whole posterior distribution on $k$. It can also be slightly modified into the smallest integer $k$ such that the *Bayes factor* comparing $\Theta_{k+1}$ to $\Theta_k$ is less than one. When considering a model comparison point of view, Bayes factors are often used to compare two models; see [11]. In the following, we shall focus on $\widehat{k}_n^{\mathrm{G}}$ and $\widehat{k}_n^{\mathrm{L}}$, since the sequential Bayes factor estimator shares the same properties as $\widehat{k}_n^{\mathrm{L}}$.

1.2. *Results in perspective.* In this paper we prove that the underestimation errors are $O(e^{-an})$ (some $a > 0$); see Theorem 1. We also show that the overestimation errors are $O((\log n)^b/n^c)$ (some $b \geq 0$, $c > 0$); see Theorems 2 and 3. All constants can be expressed explicitly, even though they are quite complicated. We apply these results in a regression model and in a change points problem. Finally, we show that our results apply to the important class of mixture models. Mixture models have interesting nonregularity properties and, in particular, even though the mixing distribution is identifiable, testing on the order of the model has proved to be difficult; see, for instance, [6]. There, we obtain an underestimation error of order $O(e^{-an})$ and an overestimation error of order $O((\log n)^b/\sqrt{n})$ ($b > 0$); see Theorem 4.

Efficiency issues in the order estimation problem have been studied mainly in the frequentist literature; see [4] for a review on these results. There is an extensive literature on Bayesian estimation of mixture models and, in particular, on the order selection in mixture models. However, this literature is essentially devoted to determining coherent noninformative priors (see, e.g., [15]) and to implementation (see, e.g., [14]). To the best of our knowledge, there is hardly any work on frequentist properties of Bayesian estimators such as $\widehat{k}_n^{\mathrm{G}}$ and $\widehat{k}_n^{\mathrm{L}}$ outside the regular case. In the case of mixture models, Ishwaran, James and San [10] suggest a Bayesian estimator of the mixing distribution when the number of components is unknown and bounded and study the asymptotic properties of the mixing distribution. It is to be noted that deriving rates of convergence for the order of the model from those of the mixing distribution would be suboptimal since the mixing distribution converges at a rate at most equal to $n^{-1/4}$ to be compared to our $O((\log n)^b/\sqrt{n})$ $(b>0)$ in Theorem 4.

1.3. *Organization of the paper.* In Section 2 we state our main results. General bounds are presented in Sections 2.1 (underestimation) and 2.2 (overestimation). The regression and change points examples are treated in Section 2.3. We deal with mixture models in Section 2.4. The main proofs are gathered in Section 3 (underestimation), Section 4 (overestimation) and Section 5 (examples). Section C in the Appendix is devoted to an aspect of mixture models which might be of interest in its own.

## 2. Efficiency bounds.

Hereafter, the integral $\int f\,d\lambda$ of a function $f$ with respect to a measure $\lambda$ is written as $\lambda f$.

Let $L^1_+(\mu)$ be the subset of all nonnegative functions in $L^1(\mu)$. For every $f\in L^1_+(\mu)\setminus\{0\}$, the measure $P_f$ is defined by its derivative $f$ with respect to $\mu$. For every $f,f'\in L^1_+(\mu)$, we set $V(f,f')=P_f(\log f-\log f')^2$ [with convention $V(f,f')=\infty$ whenever necessary].

Let $\ell^\star=\log f^\star$. For all $\theta,\theta'\in\Theta_\infty$, we set $\ell_\theta=\log f_\theta$ and define $H(\theta,\theta')=P_\theta(\ell_\theta-\ell_{\theta'})$ when $P_\theta\ll P_{\theta'}$ ($\infty$ otherwise), the Kullback–Leibler divergence between $P_\theta$ and $P_{\theta'}$. We also set $H(\theta)=H(\theta^\star,\theta)$ (each $\theta\in\Theta_\infty$).

Let us define, for every $k\geq 1$, $\alpha,\delta>0$ and $t\in\Theta_k$, $\theta\in\Theta_\infty$,

$$l_{t,\delta}=\inf\{f_{\theta'}:\theta'\in\Theta_k, d(t,\theta')<\delta\},\qquad u_{t,\delta}=\sup\{f_{\theta'}:\theta'\in\Theta_k, d(t,\theta')<\delta\},$$

$$H_k^\star=\inf\{H(\theta'):\theta'\in\Theta_k\},\qquad S_k(\delta)=\{\theta'\in\Theta_k:H(\theta')\leq H_k^\star+\delta/2\},$$

$$q(\theta,\alpha)=P^\star(\ell^\star-\ell_\theta)^2e^{\alpha(\ell^\star-\ell_\theta)}+V(\ell^\star,\ell_\theta)\in[0,\infty].$$

Throughout this paper we suppose that the following standard conditions are satisfied: for every $k\geq 1$, $(\Theta_k,d)$ is compact and $\theta\mapsto\ell_\theta(z)$ from $\Theta_k$ to $\mathbb{R}$ is continuous for every $z\in\mathcal{Z}$. By definition of $k^\star$, we have $H_k^\star=0$ for all $k\geq k^\star$ and $H_k^\star>0$ otherwise.

We consider now two assumptions that are useful for controlling the underestimation and overestimation errors.

**A1.** For each $k \geq 1$, there exist $\alpha, \delta_0 > 0$, $M \geq 1$ such that, for all $\delta \in (0, \delta_0]$,

$$\sup\{q(\theta, \alpha) : \theta \in S_k(\delta)\} \leq M.$$

**A2.** For every $k \geq 1$ and $\theta \in \Theta_k$, there exists $\eta_\theta > 0$ such that

$$V(u_{\theta,\eta_\theta}, f^\star) + V(f^\star, l_{\theta,\eta_\theta}) + V(f^\star, u_{\theta,\eta_\theta}) + V(u_{\theta,\eta_\theta}, f_\theta) < \infty.$$

Assumption **A1** states the existence of *some* (rather than *any*) exponential moment for log ratios of densities ($\ell^\star - \ell_\theta$) for $\theta$ ranging over some neighborhood of $\theta^\star$ and was also considered in [4].

2.1. *Underestimation.* We first deal with the underestimation errors.

THEOREM 1. *Assume that* **A1** *and* **A2** *are satisfied and that* $\pi_k\{S_k(\delta)\} > 0$ *for all* $\delta > 0$ *and* $k = 1, \ldots, k^\star$.

(i) *There exist* $c_1', c_2' > 0$ *such that, for every* $n \geq 1$,

$$P^{\star n}\{\widehat{k}_n^{\mathrm{G}} < k^\star\} \leq c_1' e^{-nc_2'}. \tag{1}$$

(ii) *If, in addition,* $H_k^\star > H_{k+1}^\star$ *for* $k = 1, \ldots, k^\star - 1$, *then there exist* $c_1, c_2 > 0$ *such that, for every* $n \geq 1$,

$$P^{\star n}\{\widehat{k}_n^{\mathrm{L}} < k^\star\} \leq c_1 e^{-nc_2}. \tag{2}$$

The proof of Theorem 1 is postponed to Section 3.

According to (1) and (2), both underestimation probabilities decay exponentially quickly. This is the best achievable rate. This comes from a variant of the Stein lemma (see Theorem 2.1 in [2] and Lemma 3 in [4]).

Values of constants $c_1, c_1', c_2, c_2'$ can be found in the proof of Theorem 1. Evaluating them is difficult [see (9) for a lower bound on $c_2$ in the regression model]. However, we think that they shed some light on the underestimation phenomenon. It is natural to compare our underestimation exponents $c_2$ and $c_2'$ to the constant that appears in Stein's lemma, namely, $\inf_{\theta \in \Theta_{k^\star - 1}} H(\theta, \theta^\star)$. The constants do not match, which does not necessarily mean that $\widehat{k}_n^{\mathrm{G}}$ and $\widehat{k}_n^{\mathrm{L}}$ are not optimal. We refer to [4] for a discussion about optimality.

2.2. *Overestimation.* Let the largest integer which is strictly smaller than $a \in \mathbb{R}$ be denoted by $\lfloor a \rfloor$. For simplicity, let $a \vee b$ and $a \wedge b$ be the maximum and minimum of $a, b \in \mathbb{R}$, and $V(\theta) = V(f^\star, f_\theta) \vee V(f_\theta, f^\star)$ ($\theta \in \Theta_\infty$). It is crucial in our study of overestimation errors that, if **A1** is satisfied and $C_1 = 5(1 + \log^2 M)/2\alpha^2$, then (following Lemma 5 and Theorem 5 in [20]) for all $k \geq k^\star$ and $\theta \in \Theta_k$, $H(\theta) \leq e^{-2}$ yields

$$V(\theta) \leq C_1 H(\theta) \log^2 H(\theta). \tag{3}$$

Let us now introduce further notions and assumptions. Given $\delta > 0$ and two functions $l \le u$, the bracket $[l, u]$ is the set of all functions $f$ with $l \le f \le u$. We say that $[l, u]$ is a $\delta$-bracket if $l, u \in L^1_+(\mu)$ and

$$\mu(u-l) \le \delta, \qquad P^\star(\log u - \log l)^2 \le \delta^2,$$
$$P_{u-l}(\log u - \log f^\star)^2 \le \delta \log^2 \delta \quad \text{and} \quad P_l(\log u - \log l)^2 \le \delta \log^2 \delta.$$

For $\mathcal{C}$ a class of functions, the $\delta$-entropy with bracketing of $\mathcal{C}$ is the logarithm $\mathcal{E}(\mathcal{C}, \delta)$ of the minimum number of $\delta$-brackets needed to cover $\mathcal{C}$. A set of cardinality $\exp(\mathcal{E}(\mathcal{C}, \delta))$ of $\delta$-brackets which covers $\mathcal{C}$ is written as $\mathcal{H}(\mathcal{C}, \delta)$.

For all $\theta \in \Theta_\infty$, we introduce the following quantities: $\ell_n(\theta) = \sum_{i=1}^n \ell_\theta(Z_i)$, $\ell_n^\star = \sum_{i=1}^n \ell^\star(Z_i)$ and, for every $k \ge 1$, $\mathbb{B}_n(k) = \pi(k) \int_{\Theta_k} e^{\ell_n(\theta) - \ell_n^\star} \, d\pi_k(\theta)$. Obviously, if $k < k'$ are two integers, then $\widehat{k}_n^{\mathrm{L}} = k$ yields $\mathbb{B}_n(k) \ge \mathbb{B}_n(k+1)$ and $\widehat{k}_n^{\mathrm{G}} = k$ implies that $\mathbb{B}_n(k) \ge \mathbb{B}_n(k')$.

Let $K > k^\star$ be an integer. We consider the following three assumptions:

**O1**($K$). There exist $C_2, D_1(k) > 0$ ($k = k^\star + 1, \ldots, K$) such that, for every sequence $\{\delta_n\}$ decreasing to 0, for all $n \ge 1$, and all $k \in \{k^\star + 1, \ldots, K\}$,

$$\pi_k\{\theta \in \Theta_k : H(\theta) \le \delta_n\} \le C_2 {\delta_n}^{D_1(k)/2}.$$

**O2**($K$). There exists $C_3 > 0$ such that, for each $k \in \{k^\star + 1, \ldots, K\}$, there exists a sequence $\{\mathcal{F}_n^k\}$, $\mathcal{F}_n^k \subset \Theta_k$, such that, for all $n \ge 1$,

$$\pi_k\{(\mathcal{F}_n^k)^c\} \le C_3 n^{-D_1(k)/2}.$$

**O3**. There exist $\beta_1, L, D_2(k^\star) > 0$, and $\beta_2 \ge 0$ such that, for all $n \ge 1$,

$$P^{\star n}\{\mathbb{B}_n(k^\star) < (\beta_1 (\log n)^{\beta_2} n^{D_2(k^\star)/2})^{-1}\} \le L \frac{(\log n)^{3D_1(k^\star+1)/2 + \beta_2}}{n^{[D_1(k^\star+1) - D_2(k^\star)]/2}}.$$

When **O3** holds, let $n_0$ be the smallest integer $n$ such that

$$\delta_0 = 4 \max_{m \ge n} \{m^{-1} \log[\beta_1 (\log m)^{\beta_2} m^{D_2(k^\star)/2}]\} \le e^{-2}/2. \tag{4}$$

When **O1(**$K$**)** and **O3** hold with $D_2(k^\star) < \min_{k^\star < k \le K} D_1(k)$, given any $s > 0$, we set $\delta_{k,n} = \delta_{k,1} n^{-1} \log^3 n$ for all $n \ge 2$, $k \in \{k^\star + 1, \ldots, K\}$, with

$$\delta_{k,1} \ge 128(1+s)(C_1 + 2)[D_1(k) - D_2(k^\star)] \vee 128 C_1 D_1(k) \vee \log^{-3} n_0. \tag{5}$$

We control the overestimation error for $\widehat{k}_n^{\mathrm{G}}$ when a prior bound $k_{\max}$ on $k^\star$ is known.

Theorem 2. *If the prior $\Pi$ puts mass 1 on $\bigcup_{k \le k_{\max}} \Theta_k$ and if $k^\star \le k_{\max}$, if* **A1**, **A2**, **O1(**$k_{\max}$**)**, **O2(**$k_{\max}$**)** *and* **O3** *are satisfied with $D_2(k^\star) <$* $\min_{k^\star<k\le k_{\max}} D_1(k)$, *if, in addition, for every* $k\in\{k^\star+1,\dots,k_{\max}\}$, *for all integers* $n\ge n_0$ *such that* $\delta_{k,n}<\delta_0$ *and for every* $j\le\lfloor\delta_0/\delta_{k,n}\rfloor$,

$$\mathcal{E}\left(\mathcal{F}_n^k\cap[S_k(2(j+1)\delta_{k,n})\setminus S_k(2j\delta_{k,n})],\frac{j\delta_{k,n}}{4}\right)\le\frac{s/(1+s)nj\delta_{k,n}}{64(C_1+2)\log^2(j\delta_{k,n})},\tag{6}$$

*then there exists* $c_3'>0$ *such that, for all* $n\ge n_0$,

$$P^{\star n}\{\widehat{k}_n^{\mathrm{G}}>k^\star\}\le c_3'\frac{(\log n)^{3\max_k D_1(k)/2+\beta_2}}{n^{\min_k[D_1(k)-D_2(k^\star)]/2}}.\tag{7}$$

*In the formula above index* $k$ *ranges between* $k^\star+1$ *and* $k_{\max}$.

On the contrary, the following result on the overestimation error of $\widehat{k}_n^{\mathrm{L}}$ does not rely on a prior bound on $k^\star$.

THEOREM 3. *Let* $k=k^\star+1$. *Let us suppose that assumptions* **A1**, **A2**, **O1**$(k)$, **O2**$(k)$ *and* **O3** *are satisfied with* $D_2(k^\star)<D_1(k)$. *If, in addition, for all integers* $n\ge n_0$ *such that* $\delta_{k,n}<\delta_0$ *and for every* $j\le\lfloor\delta_0/\delta_{k,n}\rfloor$, *equation (6) is satisfied, then there exists* $c_3>0$ *such that, for all* $n\ge n_0$,

$$P^{\star n}\{\widehat{k}_n^{\mathrm{L}}>k^\star\}\le c_3\frac{(\log n)^{3D_1(k^\star+1)/2+\beta_2}}{n^{[D_1(k^\star+1)-D_2(k^\star)]/2}}.\tag{8}$$

Proofs of Theorems 2 and 3 rely on tests of $P^\star$ versus complements $\{P_\theta:\theta\in\Theta_k,H(\theta)\ge\varepsilon\}$ of Kullback–Leibler balls around $P^\star$ for $k>k^\star$, in the spirit of [8]. They are postponed to Section 4. The upper bounds we get in the proofs are actually tighter than the one stated in the theorems. Each time, we actually chose the largest of several terms to make the formulas more readable. Besides, the possibility in Theorem 3 to tune the value of $\delta_{k,1}$ makes it easier to apply the theorem to the mixture model example. Naturally, the larger $\delta_{k,1}$, the larger $c_3$ and the less accurate the overestimation bound.

Concerning condition (6), it warrants that (a critical region of) $\Theta_k$ is not too large, since the entropy is known to quantify the complexity of a model.

Assumption **O1** is concerned with the decay to 0 of the prior mass of shrinking Kullback–Leibler neighborhoods of $\theta^\star$. Verifying this assumption in the mixture setting is a demanding task; see Section 2.4. Note that dimensional indices $D_1(k)$ $(k>k^\star)$ are introduced, which might be different from the usual dimensions $D(k)$. They should be understood as *effective* dimensions of $\Theta_k$ relative to $\Theta_{k^\star}$. In models of mixtures of $g_\gamma$ densities ($\gamma\in\Gamma\subset\mathbb{R}^d$), for instance, $D_1(k^\star+1)=D(k^\star)+1$, while $D(k^\star+1)=D(k^\star)+(d+1)$. It is to be noted that this assumption is crucial. In particular, in the different context of [16], it is proved that if such a condition is not satisfied, then some inconsistency occurs for the Bayes factor.

Finally, **O3** is milder than the existence of a Laplace expansion of the marginal likelihood (which holds in "regular models" as described in [18]), since in such cases (see [18]), for $c$ as large as need be, denoting by $J_n$ the Jacobian matrix, there exist $\delta, C > 0$ such that

$$\mathbb{B}_n(k^\star) \geq \int_{|\theta-\widehat{\theta}|_1 \leq \delta} e^{\ell_n(\theta)-\ell_n(\widehat{\theta})}\, d\pi_{k^\star}(\theta) \geq \left(\frac{2\pi}{n}\right)^{D(k^\star)/2} |J_n|^{-1/2}(1 + O_P(1/n)),$$

and $P^{\star n}\{|J_n| + |O_P(1/n)| > C\} \leq n^{-c}$, implying **O3** with $\beta_1 > 0$, $\beta_2 = 0$ and $D_2(k^\star) = D(k^\star)$. In some cases however, dimensional index $D_2(k^\star)$ may differ from $D(k^\star)$; see, for instance, Lemma 1.

According to (7) and (8), both overestimation errors decay as a negative power of the sample size $n$ (up to a power of a $\log n$ factor). Note that the overestimation rate is necessarily slower than exponential, as stated in another variant of the Stein Lemma (see Lemma 3 in [4]).

We want to emphasize that the overestimation rates obtained in Theorems 2 and 3 depend on intrinsic quantities [such as dimensions $D_1(k)$ and $D_2(k^\star)$, power $\beta_2$]. On the contrary, the rates obtained in Theorems 10 and 11 of [4] depend directly on the choice of a penalty term.

2.3. *Regression and change points models.* Theorems 1, 2 and 3 (resp. 1 and 3) apply to the following regression (resp. change points) model. In the rest of this section, $\sigma > 0$ is given, $g_\gamma$ is the density of the Gaussian distribution with mean $\gamma$ and variance $\sigma^2$; $X_1, \ldots, X_n$ are i.i.d. and uniformly distributed on $[0,1]$, $e_1, \ldots, e_n$ are i.i.d. with density $g_0$ and independent from $X_1, \ldots, X_n$. Moreover, one observes $Z_i = (X_i, Y_i)$ with $Y_i = \varphi_{\theta^\star}(X_i) + e_i$ $(i = 1, \ldots, n)$, where the definition of $\varphi_{\theta^\star}$ depends on the example.

*Regression* (see also Section 5.3 of [4]). Let $\{t_k\}_{k\geq 1}$ be a uniformly bounded system of continuous functions on $[0,1]$ forming an orthonormal system in $L^2([0,1])$ (for the Lebesgue measure). Let $\Gamma$ be a compact subset of $\mathbb{R}$ that contains 0 and $\Theta_k = \Gamma^k$ (each $k \geq 1$). For every $\theta \in \Theta_k$, set $\varphi_\theta = \sum_{j=1}^k \theta_j t_j$ and $f_\theta(z) = g_{\varphi_\theta(x)}(y)$ [all $z = (x,y) \in [0,1] \times \mathbb{R}$].

*Change points.* For each $k \geq 1$, let $\mathcal{T}_k$ be the set of $(k+1)$-tuples $(t_j)_{0\leq j\leq k}$, with $t_0 = 0$, $t_j \leq t_{j+1}$ (all $j < k$), and $t_k = 1$. Let $\Gamma$ be a compact subset of $\mathbb{R}$ and $\Theta_k = \mathcal{T}_k \times \Gamma^k$ (each $k \geq 1$). For every $\theta = (\alpha, t) \in \Theta_k$, set $\varphi_\theta(x) = \sum_{j=1}^k \alpha_j \mathbb{1}\{t_{j-1} \leq x < t_j\}$, and $f_\theta(z) = g_{\varphi_\theta(x)}(y)$ (all $z = (x,y) \in [0,1] \times \mathbb{R}$).

*In both examples* there exists $\theta^\star \in \Theta_{k^\star} \setminus \Theta_{k^\star - 1}$ such that $f^\star = f_{\theta^\star}$. The standard conditions of compactness and continuous parameterization are fulfilled, and **A1** and **A2** are satisfied. Besides, $2\sigma^2 H(\theta) = \|\varphi_\theta - \varphi^\star\|_2^2$ (all $\theta \in \Theta_\infty$), so the additional condition stated in Theorem 1(ii) holds. Consequently, if $\pi_k$ is positive on $\Theta_k$ for each $k \geq 1$, then Theorem 1 applies. In particular, using Fourier basis in the regression model, we get

$$12 c_2 \geq 1 / \max_{k<k^\star} \left( \frac{1}{2\sigma^2} + \frac{\Delta_{k+1}}{2\sigma^2} + \frac{2^{k^\star}}{\pi} (1 + \Delta_{k+1})^2 \right), \tag{9}$$

where $\Delta_{k+1} = (\theta^\star_{k+1})^{-2} \sum_{j=k+2}^{k^\star} (\theta^\star_j)^2$ if $k+1 < k^\star$ and $\Delta_{k^\star} = 0$.

Also, it can be shown that there exists $\tau \geq 1$ such that $[l_{\theta,\delta/\tau}; u_{\theta,\delta/\tau}]$ is a $\delta$-bracket for all $\theta \in \Theta_\infty$ and $\delta$ sufficiently small. Consequently, with the notation of Theorems 2 and 3, and with $\mathcal{F}_n^k = \Theta_k$ (**O2** is then trivial), $\mathcal{E}(\Theta_k, j\delta_{k,n}/4) \leq -b \log(j\delta_{k,n}) + c$ for positive $b, c$, and we show in Appendix D how this implies the desired condition on entropy.

*The regression model* is regular (as described in [18]), so **O3** holds with $D_2(k^\star) = D(k^\star)$. Moreover, the form of $H(\theta)$ makes it easy to verify that $\mathbf{O1}(K)$ is satisfied for any $K > k^\star$ with $D_1 = D$. Thus, Theorems 2 and 3 apply too. Furthermore, Theorem 3 applies in the *change points model* because, for any $\tau \in (0, \frac{1}{2})$ (see Appendix A for a proof),

LEMMA 1. *In the change points model,* $\mathbf{O1}(k^\star + 1)$ *and* **O3** *hold with* $D_1(k^\star + 1) = D(k^\star) + k^\star$, $D_2(k^\star) = D(k^\star) + k^\star - 1 + 2\tau$ *and* $\beta_2 = 0$.

Actually, the proof of Lemma 1 can easily be adapted to yield that $\mathbf{O1}(K)$ holds for any $K > k^\star$ with $D_1(K) = D(k^\star) + K - 1$ (we omit the details for the sake of conciseness). So Theorem 2 also applies in that model.

2.4. *Mixture models.* We prove that Theorems 1 and 3 apply here with $D_1(k^\star + 1) = D(k^\star) + 1$ and $D_2(k^\star) = D(k^\star)$, yielding an overestimation rate of order $O((\log n)^c/\sqrt{n})$ for some positive $c$.

We denote by $|\cdot|_1$ and $|\cdot|_2$ the $\ell^1$ and $\ell^2$ norms on $\mathbb{R}^d$. Let $\Gamma$ be a compact subset of $\mathbb{R}^d$ and $\Delta = \{\mathbf{g} = (g_1, \ldots, g_k) \in \Gamma^k : \min_{j<j'} |g_j - g_{j'}|_2 = 0\}$. For all $\gamma \in \Gamma$, let $g_\gamma$ be a density. In this section mixtures of $g_\gamma$'s are studied. Formally, $\Theta_1 = \Gamma$ and for every $k \geq 2$,

$$\Theta_k = \left\{\theta = (\boldsymbol{p}, \boldsymbol{\gamma}) : \boldsymbol{p} = (p_1, \ldots, p_{k-1}) \in \mathbb{R}_+^{k-1}, \sum_{j=1}^{k-1} p_j \leq 1, \boldsymbol{\gamma} \in \Gamma^k\right\}.$$

Every $\theta \in \Theta_k$ $(k \geq 2)$ is associated with $f_\theta = \sum_{j=1}^{k-1} p_j g_{\gamma_j} + (1 - \sum_{j=1}^{k-1} p_j) g_{\gamma_k}$.

Note that $D(k) = k(d+1) - 1$ for each $k \geq 1$. Also, the standard conditions of compactness and continuous parameterization are fulfilled.

We consider the following six assumptions which will be used in the mixture case. The first-, second- and third-order differentiation (with respect to $\gamma$) operators are denoted by $\nabla$, $D^2$ and $D^3$, and $|\cdot|$ stands for any norm on the space of second and third-order derivatives. We say that a function is $\mathcal{C}^k$ if it is $k$ times continuously differentiable:

**M1.** For each $k \geq 1$, prior $\pi_k$ writes as $d\pi_k(\theta) = \pi_k^{\boldsymbol{p}}(\boldsymbol{p}) \pi_k^{\boldsymbol{\gamma}}(\boldsymbol{\gamma})\, d\boldsymbol{p}\, d\boldsymbol{\gamma}$ [all $\theta = (\boldsymbol{p}, \boldsymbol{\gamma}) \in \Theta_k$]. It is $\mathcal{C}^1$ over $\Theta_k$. Moreover, there exist $\varepsilon, C > 0$ such that, setting $\Delta_\varepsilon = \{\boldsymbol{\gamma} \in \Gamma : \inf_{\mathbf{g} \in \Delta} |\boldsymbol{\gamma} - \mathbf{g}|_1 \geq \varepsilon\}$, $\boldsymbol{\gamma} \in \Delta_\varepsilon$ yields $\pi_k^{\boldsymbol{\gamma}}(\boldsymbol{\gamma}) \geq C$, and *when* $d = 1$, $\pi_k^{\boldsymbol{\gamma}}(\boldsymbol{\gamma}) \propto \prod_{j<j'} |\gamma_j - \gamma_{j'}|_2$ upon $\Delta_\varepsilon$.

**M2.** For all $\gamma \in \Gamma, \eta > 0$, let us define $\underline{g}_{\gamma,\eta} = \inf\{g_{\gamma'} : |\gamma - \gamma'|_1 \leq \eta\}$ and $\overline{g}_{\gamma,\eta} = \sup\{g_{\gamma'} : |\gamma - \gamma'|_1 \leq \eta\}$. There exist $\eta_1, M > 0$ such that, for every $\gamma_1, \gamma_2 \in \Gamma$, there exists $\eta_2 > 0$ such that

$$P_{\overline{g}_{\gamma_1,\eta_1} - \underline{g}_{\gamma_1,\eta_1}}(1 + \log^2 g_{\gamma_2}) \leq M\eta_1, \qquad P_{g_{\gamma_2}} \log^2(\overline{g}_{\gamma_1,\eta_1}/\underline{g}_{\gamma_1,\eta_1}) \leq M\eta_1^2,$$

$$P_{\overline{g}_{\gamma_1,\eta_1} - \underline{g}_{\gamma_1,\eta_1}} \log^2 \overline{g}_{\gamma_1,\eta_1} \leq M\eta_1 \log^2 \eta_1$$

and

$$P_{\overline{g}_{\gamma_1,\eta_1}}(\log^2 \overline{g}_{\gamma_1,\eta_1} + \log^2 g_{\gamma_2}) + P_{g_{\gamma_2}}(\log^2 \overline{g}_{\gamma_1,\eta_1} + \log^2 \underline{g}_{\gamma_1,\eta_1}) \leq M.$$

**M3.** For every $\gamma_1, \gamma_2 \in \Gamma$, there exists $\alpha > 0$ such that

$$\sup\{P_{\gamma_1}(g_{\gamma_2}/g_\gamma)^\alpha : \gamma \in \Gamma\} < \infty.$$

**M4.** The parameterization $\gamma \mapsto g_\gamma(z)$ is $\mathcal{C}^2$ for $\mu$-almost every $z \in \mathcal{Z}$. Moreover, $\mu[\sup_{\gamma\in\Gamma}(|\nabla g_\gamma|_1 + |D^2 g_\gamma|)]$ is finite.

The parameterization $\gamma \mapsto \log g_\gamma(z)$ is $\mathcal{C}^3$ for $\mu$-almost every $z \in \mathcal{Z}$ and for every $\gamma \in \Gamma$, the Fisher information matrix $I(\gamma)$ is positive definite. Besides, for all $\gamma_1, \gamma_2 \in \Gamma$, there exists $\eta > 0$ for which

$$P_{\gamma_1}|D^2 \log g_{\gamma_2}|^2 + P_{\gamma_1} \sup\{|D^3 \log g_\gamma|^2 : |\gamma - \gamma_2|_1 \leq \eta\} < \infty.$$

**M5.** Let $\mathcal{I} = \{(r,s) : 1 \leq r \leq s \leq d\}$. There exist a nonempty subset $\mathcal{A}$ of $\mathcal{I}$ and two constants $\eta_0, a > 0$ such that, for every $k \geq 2$, for every $k$-tuple $(\gamma_1, \ldots, \gamma_k)$ of pairwise distinct elements of $\Gamma$:

(a) functions $g_{\gamma_j}, (\nabla g_{\gamma_j})_l$ $(j \leq k, l \leq d)$ are linearly independent;

(b) for every $j \leq k$, functions $g_{\gamma_j}, (\nabla g_{\gamma_j})_l, (D^2 g_{\gamma_j})_{rs}$ (all $l \leq d, (r,s) \in \mathcal{A}$) are linearly independent;

(c) for each $j \leq k, (r,s) \in \mathcal{I} \setminus \mathcal{A}$, there exist $\lambda_{rs}^{0j}, \ldots, \lambda_{rs}^{dj} \in \mathbb{R}$ such that $(D^2 g_{\gamma_j})_{rs} = \lambda_{rs}^{0j} g_{\gamma_j} + \sum_{l=1}^d \lambda_{rs}^{lj} (\nabla g_{\gamma_j})_l$;

(d) for all $\eta \leq \eta_0$ and all $u, v \in \mathbb{R}^d$, for each $j \leq k$, if

$$\sum_{(r,s)\in\mathcal{A}} (|u_r u_s| + |v_r v_s|) + \left| \sum_{(r,s)\notin\mathcal{A}} \lambda_{rs}^{0j}(u_r u_s + v_r v_s) \right| \leq \eta,$$

then $|u|_2^2 + |v|_2^2 \leq a\eta$.

These assumptions suffice to guarantee the bounds below.

THEOREM 4. *If* **M1–M5** *are satisfied, then there exists* $n_1 \geq 1$ *and* $c_4 > 0$ *such that, for all* $n \geq n_1$,

$$P^\star\{\widehat{k}_n^{\mathrm{L}} < k^\star\} \leq c_1 e^{-nc_2}, \tag{10}$$

$$P^\star\{\widehat{k}_n^{\mathrm{L}} > k^\star\} \leq c_4 \frac{(\log n)^{[3(d+1)k^\star/2]}}{\sqrt{n}}. \tag{11}$$

*The positive constants* $c_1, c_2$ *are defined in Theorem 1.*

Note that all assumptions involve the mixed densities $g_\gamma$ ($\gamma \in \Gamma$) rather than the resulting mixture densities $f_\theta$ ($\theta \in \Theta_\infty$). Assumption **M2** implies **A2** and **M3** implies **A1**. Assumption **M4** is a usual regularity condition. Assumption **M5** is a weaker version of the strong identifiability condition defined by [5], which is assumed in most paper dealing with asymptotic properties of mixtures. In particular, strong identifiability does not hold in location-scale mixtures of Gaussian r.v., but **M5** does (with $\mathcal{A} = \mathcal{I} \setminus \{(1,1)\}$). In fact, Theorem 4 applies, and we have the following:

COROLLARY 1. *Set $A, B > 0$ and $\Gamma = \{(\mu, \sigma^2) \in [-A, A] \times [\frac{1}{B}, B]\}$. For every $\gamma = (\mu, \sigma^2) \in \Gamma$, let us denote by $g_\gamma$ the Gaussian density with mean $\mu$ and variance $\sigma^2$. Then (10) and (11) hold with $d = 2$ for all $n \geq n_0$.*

Other examples include, for instance, mixtures of Gamma$(a, b)$ in $a$ or in $b$ [but not in $(a, b)$], of Beta$(a, b)$ in $(a, b)$, of GIG$(a, b, c)$ in $(b, c)$ (another example where strong identifiability does not hold, but **M5** does).

**3. Underestimation proofs.** Let us start with new notation. For $f, f' \in L^1_+(\mu) \setminus \{0\}$, we set $H(f, f') = P_f(\log f - \log f')$ when it is defined ($\infty$ otherwise), $H(f) = H(f^\star, f)$, and $V(f) = V(f^\star, f) \vee V(f, f^\star)$. For every $\theta \in \Theta_\infty$, the following shortcuts will be used ($W$ stands for $H$ or $V$): $W(f, f_\theta) = W(f, \theta)$, $W(f_\theta, f) = W(\theta, f)$, $W(f_\theta) = W(\theta)$. For every probability density $f \in L^1(\mu)$, $P_f^{\otimes n}$ is denoted by $P_f^n$ and the expectation with respect to $P_f$ (resp. $P_f^n$) by $\mathrm{E}_f$ (resp. $\mathrm{E}_f^n$).

Theorem 1 relies on the following lower bound on $\mathbb{B}_n(k)$.

LEMMA 2. *Let $k \leq k^\star$ and $\delta \in (0, \alpha M \wedge \delta_0]$. Under the assumptions of Theorem 1, with probability at least $1 - 2\exp\{-n\delta^2/8M\}$,*

$$\mathbb{B}_n(k) \geq \frac{\pi(k)\pi_k\{S_k(\delta)\}}{2} e^{-n[H_k^\star + \delta]}.$$

PROOF. Let $1 \leq k \leq k^\star$, $0 < \delta \leq \alpha M \wedge \delta_0$ and define

$$B = \{(\theta, Z^n) \in \Theta_k \times \mathcal{Z}^n : \ell_n(\theta) - \ell_n^\star \geq -n[H_k^\star + \delta]\}.$$

Then, using the same calculations as in Lemma 1 of [17], we obtain

$$P^{\star n}\left\{\pi_k\{S_k(\delta) \cap B^c\} \geq \frac{\pi_k\{S_k(\delta)\}}{2}\right\} \leq \int_{S_k(\delta)} \frac{2P^{\star n}\{B^c\}}{\pi_k\{S_k(\delta)\}} d\pi_k(\theta). \tag{12}$$

Set $s \in [0, \alpha]$ and $\theta \in S_k(\delta)$ and let $\varphi_\theta(t) = P^\star e^{t(\ell^\star - \ell_\theta)}$ (every $t \in \mathbb{R}$). By virtue of **A1**, function $\varphi_\theta$ is $\mathcal{C}^\infty$ over $[0, \alpha]$ and $\varphi_\theta''$ is bounded by $q(\theta, \alpha) \leq M$ on that interval. Moreover, a Taylor expansion implies that

$$\varphi_\theta(s) = 1 + sH(\theta) + s^2 \int_0^1 (1 - t)\varphi_\theta''(st)\, dt \leq 1 + sH(\theta) + \tfrac{1}{2}s^2 M.$$

Applying the Chernoff method and inequality $\log t \le t-1$ $(t>0)$ implies that, for all $\theta \in S_k(\delta)$,

$$P^{\star n}\{B^c\} \le \exp\{-ns[H_k^\star + \delta] + n\log\varphi_\theta(s)\}$$
$$\le \exp\{-ns[H_k^\star + \delta - H(\theta)] + ns^2 M/2\}.$$

We choose $s = [H_k^\star + \delta - H(\theta)]/M \in [\delta/2, \alpha]$ so that the above probability is bounded by $\exp\{-n\delta^2/8M\}$ and Lemma 2 is proved. □

To prove Theorem 1, we construct nets of upper bounds for the $f_\theta$'s ($\theta \in \Theta_k$, $k = 1, \ldots, k^\star - 1$). Similar nets have been first introduced in a context of nonparametric Bayesian estimation in [3]. We focus on $\widehat{k}_n^{\mathrm{L}}$; the proof for $\widehat{k}_n^{\mathrm{G}}$ is a straightforward adaptation.

Proof of Theorem 1. Since $P^{\star n}\{\widehat{k}_n^{\mathrm{L}} < k^\star\} = \sum_{k=1}^{k^\star-1} P^{\star n}\{\widehat{k}_n^{\mathrm{L}} = k\}$, it is sufficient to study $P^{\star n}\{\widehat{k}_n^{\mathrm{L}} = k\}$ for $k$ between 1 and $k^\star - 1$.

Let $\delta < \alpha M \wedge \delta_0 \wedge [H_k^\star - H_{k+1}^\star]/2$, $c = \frac{1}{2}\pi(k)\pi_k\{S_k(\delta)\} \in (0,1]$ and $\varepsilon = 2\delta/[H_k^\star - H_{k+1}^\star] \in (0,1)$. Lemma 2 yields

$$P^{\star n}\{\widehat{k}_n^{\mathrm{L}} = k\} \le P^{\star n}\{\mathbb{B}_n(k) \ge \mathbb{B}_n(k+1)\}$$
$$\le 2e^{-n\delta^2/(8M)} + P^{\star n}\{\mathbb{B}_n(k) \ge ce^{-n[H_{k+1}^\star + \delta]}\}. \tag{13}$$

We now study the rightmost term of (13). Let $\theta, \theta' \in \Theta_k$. The dominated convergence theorem and **A2** ensure that there exists $\eta_\theta > 0$ such that $d(\theta,\theta') < \eta_\theta$ yields $H(u_{\theta,\eta}) \le H(\theta') \le H(u_{\theta,\eta}) + \delta$, $V(\theta^\star, u_{\theta,\eta}) \le (1+\varepsilon)V(\theta^\star, \theta')$ and $V(u_{\theta,\eta}, \theta^\star) \le (1+\varepsilon)V(\theta', \theta^\star)$. Let $\mathcal{B}(\theta, \eta_\theta) = \{\theta' \in \Theta_k : d(\theta,\theta') < \eta_\theta\}$ for all $\theta \in \Theta_k$. The collection of open sets $\{\mathcal{B}(\theta,\eta_\theta)\}_{\theta\in\Theta_k}$ covers $\Theta_k$, which is a compact set. So, there exist $\theta_1, \ldots, \theta_{N_\varepsilon} \in \Theta_k$ such that $\Theta_k = \bigcup_{j=1}^{N_\varepsilon} \mathcal{B}(\theta_j, \eta_{\theta_j})$. For $j = 1, \ldots, N_\varepsilon$, letting $u_j = u_{\theta_j, \eta_{\theta_j}}$,

$$\widetilde{T}_{kj} = \{\theta \in \Theta_k : \ell_\theta \le \log u_j, H(\theta) \le H(u_j) + \delta,$$
$$V(\theta^\star, u_j) \le (1+\varepsilon)V(\theta^\star, \theta), V(u_j, \theta^\star) \le (1+\varepsilon)V(\theta, \theta^\star)\},$$

then $T_{k1} = \widetilde{T}_{k1}$ and $T_{kj} = \widetilde{T}_{kj} \cap (\bigcup_{j'<j} \widetilde{T}_{kj'})^c$ $(j = 2, \ldots, N_\varepsilon)$. The family $\{T_{k1}, \ldots, T_{kN_\varepsilon}\}$ is a partition of $\Theta_k$.

Accordingly, with $\ell_{n,u_j} = \sum_{i=1}^n \log u_j(Z_i)$ $(j = 1, \ldots, N_\varepsilon)$, the rightmost term of (13) is smaller than

$$P^{\star n}\left\{\sum_{j=1}^{N_\varepsilon} \int_{T_{kj}} e^{\ell_n(\theta) - \ell_n^\star}\, d\pi_k(\theta) \ge ce^{-n[H_{k+1}^\star + \delta]}\right\}$$
$$\le \sum_{j=1}^{N_\varepsilon} P^{\star n}\left\{e^{\ell_{n,u_j} - \ell_n^\star} \int_{T_{kj}} e^{\ell_n(\theta) - \ell_{n,u_j}}\, d\pi_k(\theta) \ge ce^{-n[H_{k+1}^\star + \delta]}\pi_{k+1}\{T_{kj}\}\right\}$$

$$\leq \sum_{j=1}^{N_\varepsilon} P^{\star n}\{\ell_{n,u_j} - \ell_n^\star \geq -n[H_{k+1}^\star + \delta] + \log c\}$$

$$\leq \sum_{j=1}^{N_\varepsilon} P^{\star n}\{\ell_{n,u_j} - \ell_n^\star + nH(u_j) \geq n\rho_j + \log c\}$$

for $\rho_j = [H(u_j) - H_{k+1}^\star - \delta]$. Note that $\rho_j \geq (1-\varepsilon)[H(\theta_j) - H_{k+1}^\star] > 0$ for $j = 1, \dots, N_\varepsilon$ by construction. Applying (29) of Proposition B.1 (whose assumptions are satisfied) finally implies that

$$\begin{aligned} &P^{\star n}\{\mathbb{B}_n(k) \geq ce^{-n[H_{k+1}^\star + \delta]}\} \\ &\quad \leq \frac{N_\varepsilon}{c} \exp\left\{-n\frac{(1-\varepsilon)^2}{2(1+\varepsilon)}[H_k^\star - H_{k+1}^\star] \min\left(\inf_{\theta\in\Theta_k} \frac{H(\theta) - H_{k+1}^\star}{V(\theta)}, \frac{1+\varepsilon}{1-\varepsilon}\right)\right\}. \end{aligned}$$

We conclude, since $N_\varepsilon$ does not depend on $n$. □

**4. Overestimation proofs.** We choose again to focus on $\widehat{k}_n^{\mathrm{L}}$, the proof for $\widehat{k}_n^{\mathrm{G}}$ being very similar.

PROOF OF THEOREM 3. Set $n_0$ and $\delta_0$ as in (4), then note that $u \mapsto u\log^2 u$ increases on interval $(0, e^{-2})$. By definition of $\widehat{k}_n^{\mathrm{L}}$,

$$\begin{aligned} P^{\star n}\{\widehat{k}_n^{\mathrm{L}} > k^\star\} &\leq P^{\star n}\{\mathbb{B}_n(k^\star) < \mathbb{B}_n(k^\star+1)\} \\ &\leq P^{\star n}\{\mathbb{B}_n(k^\star) \leq (\beta_1(\log n)^{\beta_2} n^{D_2(k^\star)/2})^{-1}\} \\ &\quad + P^{\star n}\{\mathbb{B}_n(k^\star+1) \geq (\beta_1(\log n)^{\beta_2} n^{D_2(k^\star)/2})^{-1}\}. \end{aligned} \tag{14}$$

Assumption **O3** deals with the first term of the right-hand side of (14). Let us focus on the second one. To this end, $\Theta_{k^\star+1}$ is decomposed into the following three sets: letting $\delta_1$ satisfy (5) and $\delta_n = \delta_1 n^{-1}\log^3 n$,

$$\begin{aligned} S_{k^\star+1}(2\delta_0)^c &= \{\theta \in \Theta_{k^\star+1} : H(\theta) > \delta_0\}, \\ S_n &= S_{k^\star+1}(2\delta_0) \cap S_{k^\star+1}(2\delta_n)^c = \{\theta \in \Theta_{k^\star+1} : \delta_n < H(\theta) \leq \delta_0\}, \\ S_{k^\star+1}(2\delta_n) &= \{\theta \in \Theta_{k^\star+1} : H(\theta) \leq \delta_n\}. \end{aligned}$$

Note that $S_n$ can be empty. According to this decomposition, the quantity of interest is bounded by the sum of three terms (the second one is 0 when $S_n$ is empty): if $w_n = 3\pi(k^\star+1)\beta_1(\log n)^{\beta_2} n^{D_2(k^\star)/2}$, then

$$\begin{aligned} &P^{\star n}\{\mathbb{B}_n(k^\star+1) \geq (\beta_1(\log n)^{\beta_2} n^{D_2(k^\star)/2})^{-1}\} \\ &\quad \leq P^{\star n}\left\{\int_{S_{k^\star+1}(2\delta_0)^c} e^{\ell_n(\theta) - \ell_n^\star}\, d\pi_{k^\star+1}(\theta) \geq 1/w_n\right\} \end{aligned}$$

(15)

$$+ P^{\star n}\Big\{\int_{S_n} e^{\ell_n(\theta)-\ell_n^\star}\, d\pi_{k^\star+1}(\theta) \geq 1/w_n\Big\}$$
$$+ P^{\star n}\Big\{\int_{S_{k^\star+1}(2\delta_n)} e^{\ell_n(\theta)-\ell_n^\star}\, d\pi_{k^\star+1}(\theta) \geq 1/w_n\Big\}.$$

The Markov inequality, Fubini theorem and **O1** yield (as in the proof of Lemma 2) the following bound on the third term, $p_{n,3}$, of (15):

$$\begin{aligned} p_{n,3} &\leq w_n \pi_{k^\star+1}\{S_{k^\star+1}(2\delta_n)\} \leq C_2 w_n {\delta_n}^{D_1(k^\star+1)/2} \\ &\leq 3\beta_1 C_2 \pi(k^\star+1){\delta_1}^{D_1(k^\star+1)/2} \frac{(\log n)^{3D_1(k^\star+1)/2+\beta_2}}{n^{[D_1(k^\star+1)-D_2(k^\star)]/2}}. \end{aligned} \tag{16}$$

The first term of (15), $p_{n,1}$, is like $P^{\star n}\{\mathbb{B}_n(k) \geq ce^{-n[H_{k+1}^\star+\delta]}\}$, already bounded in the proof of Theorem 1. Indeed, the infima for $\theta \in S_{k^\star+1}(2\delta_0)^c$ of $H(\theta)$, $V(\theta^\star,\theta)$ and $V(\theta,\theta^\star)$ are positive and the scheme of proof of Theorem 1 also applies here: there exist $c_4, c_5 > 0$ which do not depend on $n$ and guarantee that

$$p_{n,1} \leq c_4 e^{-nc_5}. \tag{17}$$

When $\delta_n < \delta_0$, bounding the second term of (15), $p_{n,2}$, goes in four steps.

Let $\Delta_n = \lfloor \delta_0/\delta_n \rfloor$. For all $j = 1,\dots,\Delta_n$, let $S_{n,j} = \{\theta \in \mathcal{F}_n \cap S_n : j\delta_n < H(\theta) \leq (j+1)\delta_n\}$. Consider $[l_i,u_i] \in \mathcal{H}(S_{n,j}, j\delta_n/4)$, define $\overline{u}_i = u_i/\mu u_i$ and introduce the local tests

$$\phi_{i,j} = \mathbb{1}\Big\{\ell_{n,\overline{u}_i} - \ell_n^\star + nH(\overline{u}_i) \geq n\frac{j\delta_n}{2}\Big\} = \phi_{n,f,\rho,c}$$

for $f = \overline{u}_i$, $\rho = j\delta_n/2$ and $c = 1$ in the perspective of Proposition B.1.

*Step* 1. Set $\theta \in S_{n,j}$ such that $f_\theta \in [l_i,u_i]$, $g = f_\theta$ and $\rho' = \log \mu u_i$. Then $\mu g = 1$, $V(g) = V(\theta) > 0$ and $H(\overline{u}_i) - (\rho+\rho') = P^\star(\ell^\star - \log \overline{u}_i) - \log \mu u_i - \rho = P^\star(\ell^\star - \log u_i) - \rho = H(\theta) + P^\star(\ell_\theta - \log u_i) - \rho \geq H(\theta) - P^\star(\log u_i - \log l_i) - \rho \geq \frac{j\delta_n}{4} > 0$. Thus, according to (30) of Proposition B.1,

$$\mathrm{E}_\theta^n(1-\phi_{i,j}) \leq \exp\{-\frac{n[H(\overline{u}_i)-(\rho+\rho')]}{2}\Big(\frac{H(\overline{u}_i)-(\rho+\rho')}{V(\theta)} \wedge 1\Big)\Big\}.$$

Since $H(\theta) \leq (j+1)\delta_n \leq 2\delta_0 \leq e^{-2}$, then $\log^2\delta_n \geq \log^2(j\delta_n)$ and (3) yield $V(\theta) \leq C_1 H(\theta)\log^2 H(\theta) \leq C_1(j+1)\delta_n \log^2(j\delta_n)$. Consequently, $j/(j+1) \geq 1/2$ and $8C_1\log^2(j\delta_n) \geq 1$ imply

$$\mathrm{E}_\theta^n(1-\phi_{i,j}) \leq \exp\Big\{-\frac{nj\delta_n}{64C_1\log^2(j\delta_n)}\Big\}. \tag{18}$$

*Step* 2. Proposition B.1 and (29) ensure that

$$\mathrm{E}^{\star n}\phi_{i,j} \leq \exp\Big\{-\frac{nj\delta_n}{4}\Big(\frac{j\delta_n}{2V(\overline{u}_i)} \wedge 1\Big)\Big\}.$$

The point is now to bound $V(\overline{u}_i)$. Let again $\theta \in S_{n,j}$ be such that $f_\theta \in [l_i, u_i]$. Using repeatedly $(a+b)^2 \le 2(a^2+b^2)$ $(a, b \in \mathbb{R})$, the definition of a $\delta$-bracket and (3) yield

$$\begin{aligned} V(\theta^\star, \overline{u}_i) &= P^\star(\ell^\star - \log u_i + \log \mu u_i)^2 \\ &\le 2P^\star(\ell^\star - \log u_i)^2 + 2\log^2 \mu u_i \\ &\le 4P^\star(\ell^\star - \ell_\theta)^2 + 4P^\star(\ell_\theta - \log u_i)^2 + 2(\mu(u_i - l_i))^2 \\ &\le 4V(\theta) + 4P^\star(\log u_i - \log l_i)^2 + 2(\mu(u_i - l_i))^2 \\ &\le 2(2C_1+3)(j+1)\delta_n \log^2(j\delta_n), \end{aligned} \tag{19}$$

and similarly,

$$V(\overline{u}_i, \theta^\star) \le 4(C_1+2)(j+1)\delta_n \log^2(j\delta_n). \tag{20}$$

A bound on $V(\overline{u}_i)$ is derived from (19) and (20), which yields in turn

$$\mathrm{E}^{\star n}\phi_{i,j} \le \exp\left\{-\frac{nj\delta_n}{64(C_1+2)\log^2(j\delta_n)}\right\}. \tag{21}$$

*Step* 3. Now, consider the global test

$$\phi_n = \max\{\phi_{i,j} : i \le \exp\{\mathcal{E}(S_{n,j}, j\delta_n/4)\}, j \le \Delta_n\}.$$

Equation (18) implies that, for every $j \le \Delta_n$ and $\theta \in S_{n,j}$,

$$\mathrm{E}^n_\theta(1-\phi_n) \le \exp\left\{-\frac{nj\delta_n}{64C_1\log^2(j\delta_n)}\right\}. \tag{22}$$

Furthermore, if we set $\rho_n = n\delta_n/[64(1+s)(C_1+2)\log^2 \delta_n]$, then bounding $\phi_n$ by the sum of all $\phi_{i,j}$, invoking (21) and (6) yield

$$\begin{aligned} \mathrm{E}^{\star n}\phi_n &\le \sum_{j=1}^{\Delta_n} \exp\left\{\mathcal{E}(S_{n,j}, j\delta_n/4) - \frac{nj\delta_n}{64(C_1+2)\log^2(j\delta_n)}\right\} \\ &\le \sum_{j=1}^{\Delta_n} \exp\{-j\rho_n\} \le \frac{\exp\{-\rho_n\}}{1-\exp\{-\rho_n\}}. \end{aligned}$$

Since $\delta_1 \ge 128(1+s)(C_1+2)[D_1(k^\star+1) - D_2(k^\star)] \vee \log^{-3}(n_0)$, one has $\log^2 \delta_n \le 4\log^2 n$, and $\rho_n \ge \frac{1}{2}[D_1(k^\star+1) - D_2(k^\star)]\log n$. Thus, the final bound is

$$\mathrm{E}^{\star n}\phi_n \le \frac{1}{n^{[D_1(k^\star+1)-D_2(k^\star)]/2} - 1}. \tag{23}$$

*Step* 4. We now bound $p_{n,2}$:

$$p_{n,2} = \mathrm{E}^{\star n}\Big(\phi_n + (1-\phi_n)\Big)\mathbb{1}\Big\{\int_{S_n} e^{\ell_n(\theta)-\ell_n^\star}\, d\pi_{k^\star+1}(\theta) \geq 1/w_n\Big\}$$
$$\leq \mathrm{E}^{\star n}\phi_n + P^{\star n}\Big\{\int_{S_n\cap\mathcal{F}_n^c} e^{\ell_n(\theta)-\ell_n^\star}\, d\pi_{k^\star+1}(\theta) \geq 1/2w_n\Big\}$$
$$+ \mathrm{E}^{\star n}(1-\phi_n)\mathbb{1}\Big\{\int_{S_n\cap\mathcal{F}_n} e^{\ell_n(\theta)-\ell_n^\star}\, d\pi_{k^\star+1}(\theta) \geq 1/2w_n\Big\}.$$

The first term of the right-hand side is bounded according to (23). Moreover, applying the Markov inequality and Fubini theorem to the second term above, $p_{n,2,2}$, ensures that

$$p_{n,2,2} \leq 6\beta_1 C_3 \frac{(\log n)^{\beta_2}}{n^{[D_1(k^\star+1)-D_2(k^\star)]/2}}. \tag{24}$$

As for the third term, $p_{n,2,3}$, invoking again the Markov inequality and Fubini theorem, then (22), yields

$$p_{n,2,3} \leq 2w_n \sum_{j=1}^{\Delta_n} \int_{S_{n,j}} \mathrm{E}_\theta^n(1-\phi_n)\, d\pi_{k^\star+1}(\theta)$$
$$\leq 2w_n \sum_{j=1}^{\Delta_n} \exp\Big\{-\frac{nj\delta_n}{64C_1\log^2(j\delta_n)}\Big\}\pi_{k^\star+1}\{S_{n,j}\}$$
$$\leq 2w_n \exp\Big\{-\frac{n\delta_n}{64C_1\log^2\delta_n}\Big\} \leq 2w_n \exp\Big\{-\frac{\delta_1}{256C_1}\log n\Big\}$$
$$\leq 6\beta_1\pi(k^\star+1)\frac{(\log n)^{\beta_2}}{n^{[D_1(k^\star+1)-D_2(k^\star)]/2}}.$$

Combining inequalities (23), (24) and (25) yields

$$p_{n,2} \leq \frac{1}{n^{[D_1(k^\star+1)-D_2(k^\star)]/2}-1} + 6\beta_1(\pi(k^\star+1)+C_3)\frac{(\log n)^{\beta_2}}{n^{[D_1(k^\star+1)-D_2(k^\star)]/2}}.$$

Inequalities (16), (17) and the one above conclude the proof. □

**5. Mixtures proofs.** In the sequel we use the notation $\theta^\star = (\boldsymbol{p}^\star, \boldsymbol{\gamma}^\star)$, $\boldsymbol{p}^\star = (p_1^\star, \ldots, p_{k^\star-1}^\star)$ and $p_{k^\star}^\star = 1 - \sum_{j=1}^{k^\star-1} p_j^\star$. Also, if $\theta = (\boldsymbol{p}, \boldsymbol{\gamma}) \in \Theta_k$, then $1 - \sum_{j=1}^{k-1} p_j$ is denoted by $p_k$.

The standard conditions hold. Assumption **A1** is verified by proving (with usual regularity and convexity arguments) the existence of $\alpha > 0$ such that the function $\theta \mapsto P^\star e^{\alpha(\ell^\star - \ell_\theta)}$ is bounded on $\Theta_{k^\star}$. Assumption **A2** follows from **M2**. Lemma 3 in [12] guarantees that $H_k^\star > H_{k+1}^\star$ (every $k < k^\star$).

So, the underestimation error bound (10) in Theorem 4 is a consequence of Theorem 1.

The overestimation error bound (11) in Theorem 4 is a consequence of Theorem 3. Let us verify that **O1**$(k^\star+1)$, **O2**$(k^\star+1)$ and **O3** are satisfied.

PROPOSITION 1. *There exists $C_2>0$ such that, in the setting of mixture models, for every sequence $\{\delta_n\}$ that decreases to 0, for all $n\geq 1$,*

$$\pi_{k^\star+1}\{\theta\in\Theta_{k^\star+1}: H(\theta)\leq\delta_n\}\leq C_2{\delta_n}^{[D(k^\star)+1]/2}.$$

PROPOSITION 2. *If $\mathcal{F}_n^{k^\star+1}=\{(\boldsymbol{p},\boldsymbol{\gamma})\in\Theta_{k^\star+1}:\min_{j\leq k^\star+1}p_j\geq e^{-n}\}$ approximates the set $\Theta_{k^\star+1}$, then* **O2**$(k^\star+1)$ *is fulfilled. Furthermore, the entropy condition* (6) *holds as soon as $\delta_1$ is chosen large enough.*

The technical proofs of Propositions 1 and 2 are postponed to Appendix C and D, respectively. Assumption **O3** is obtained (with $\beta_2=0$) from the Laplace expansion under $P^\star$, which is regular (see also the comment after Theorem 3). Finally, Theorem 3 applies and Theorem 4 is proven.

## APPENDIX A: PROOF OF LEMMA 1

Let $\theta^\star=(\alpha^\star,t^\star)$ and $\theta\in\Theta_{k^\star+1}$ satisfy $H(\theta)\leq\delta_n$. For every $j\leq k^\star$ (resp. $j\leq k$), we denote by $\tau_j^\star$ (resp. $\tau_j$) the interval $[t_{j-1}^\star,t_j^\star[$ (resp. $[t_{j-1},t_j[$) [hence, $H(\theta)=\sum_{j\leq k^\star}\sum_{j'\leq k^\star+1}(\alpha_j^\star-\alpha_{j'})^2\mu(\tau_j^\star\cap\tau_{j'})$], and set $s(j)$ such that $\mu(\tau_j^\star\cap\tau_{s(j)})=\max_{l\leq k}\mu(\tau_j^\star\cap\tau_l)$. So, $\mu(\tau_j^\star\cap\tau_{s(j)})\geq\mu(\tau_j^\star)/k$, and $(\alpha_j^\star-\alpha_{s(j)})^2\leq c\delta_n$ for all $j\leq k^\star$. If $s(j)=s(j')$ for $j'>j$, then necessarily $j'\geq(j+2)$ and $s(j+1)=s(j)$, while $\alpha_j^\star\neq\alpha_{j+1}^\star$, so we do get $k^\star$ conditions on $\theta$. Suppose now without loss of generality that $s(j)=j$ for all $j\leq k^\star$. Then $(\alpha_k-\alpha_{k^\star}^\star)^2(1-t_{k^\star})\leq\delta_n$, another condition on $\theta$. Moreover, for all $j<k^\star$, $\mu(\tau_j^\star)-\mu(\tau_j^\star\cap\tau_j)=\mu(\tau_j^\star\cap\tau_{j-1})+\mu(\tau_j^\star\cap\tau_{j+1})$, $\mu(\tau_j)-\mu(\tau_j^\star\cap\tau_j)=\mu(\tau_j\cap\tau_{j-1}^\star)+\mu(\tau_j\cap\tau_{j+1}^\star)$ (with convention $\tau_{-1}=\tau_{-1}^\star=\varnothing$) and $\alpha_j^\star\neq\alpha_{j+1}^\star$ imply $|\mu(\tau)-\mu(\tau_j^\star\cap\tau_j)|\leq c\delta_n$ for $\tau\in\{\tau_j^\star,\tau_j\}$. So, $|(t_j^\star-t_j)-(t_{j-1}^\star-t_{j-1})|\leq 2c\delta_n$. Using successively these inequalities from $j=1$ to $j=(k^\star-1)$, we get $(k^\star-1)$ conditions on $\theta$ of the form $|t_j^\star-t_j|\leq c\delta_n$. Combining those conditions yields **O1**$(k^\star+1)$ with $D_1(k^\star+1)=D(k^\star)+k^\star$.

Let $S_n=\{t^\star+u/n:u\in\mathbb{R}_+^{k^\star+1},u_0=u_{k^\star}=0,|u|_1\leq\frac{\tau}{2}\log\log n\}\subset\mathcal{T}_{k^\star}$. For large $n$, there exists an event of probability $1-(1-\min_k|t_k^\star-t_{k-1}^\star|/2)^n$ upon which the model is regular in $\alpha$ for any fixed $t\in S_n$, hence, there exists $C>0$ (independent of $t$) such that, on that event,

$$\int_{\Gamma^{k^\star}}e^{\ell_n(\theta)-\ell_n^\star}\,d\pi_{k^\star}(\alpha|t)\geq\frac{C}{n^{k^\star/2}}e^{\ell_n(\hat\alpha_t,t)-\ell_n^\star}\geq\frac{C}{n^{k^\star/2}}e^{\ell_n(\alpha^\star,t)-\ell_n^\star},\tag{25}$$

where $\hat\alpha_t$ is the maximum likelihood estimator for fixed $t$. Denote $n_j(t)=\sum_{i=1}^n\mathbb{1}\{X_i\in[t_j^\star,t_j^\star+u_j/n[\}$ and $v^2(t)=\sigma^2\sum_{j=1}^{k^\star}(\alpha_j^\star-\alpha_{j-1}^\star)^2n_j(t)$ for any $t \in S_n$. Then $\xi(t) = \ell_n(\alpha^\star, t) - \ell_n^\star + \frac{1}{2}v^2(t)$ is, conditionally on $X_1, \ldots, X_n$, a centered Gaussian r.v. with variance $v^2(t)$. Because each $n_j(t)$ is Binomial$(n, u_j/n)$ distributed, the Chernoff method implies, for any $t \in S_n$,

$$P^{\star n}\{v^2(t) \geq \tau \log n\} = O(1/\sqrt{n}). \tag{26}$$

Moreover, since $\xi(t)$ is conditionally Gaussian, it is easily seen by using (26) that, for any $t \in S_n$, setting $B = \{Z^n : \ell_n(\alpha^\star, t) - \ell_n^\star \geq -\frac{1}{2}(v^2(t) + \tau \log n)\}$,

$$P^{\star n}\{B^c\} = O(1/\sqrt{n}), \tag{27}$$

too. Now, the same technique as in the proof of Lemma 2 yields

$$P^{\star n}\left\{\int_{S_n} e^{\ell_n(\alpha^\star, t) - \ell_n^\star} \, d\pi_{k^\star}(t) \leq n^{-k^\star + 1 - \tau}\right\} \leq \int_{S_n} \frac{2P^{\star n}\{B^c\}}{\pi_{k^\star}(S_n)} \, d\pi_{k^\star}(t) \tag{28}$$

whenever $\pi_{k^\star}\{S_n\} = c(\log\log n/n)^{k^\star - 1} \geq 2n^{-k^\star + 1}$. By combining (25), (27) and (28), we obtain that **O3** holds with $D_2(k^\star) = 3k^\star + 2(\tau - 1)$.

## APPENDIX B: CONSTRUCTION OF TESTS

PROPOSITION B.1. *Let $(\rho, c)$ belong to $\mathbb{R}_+^* \times (0, 1]$ and $f \in L_+^1(\mu) \setminus \{0\}$. Assume that $V(f)$ is positive and finite. Let $\ell_{n,f} = \sum_{i=1}^n \log f(Z_i)$ and*

$$\phi_{n,f,\rho,c} = \mathbb{1}\{\ell_{n,f} - \ell_n^\star + nH(f) \geq n\rho + \log c\}.$$

*The following bound holds:*

$$\mathrm{E}^{\star n} \phi_{n,f,\rho,c} \leq \frac{1}{c} \exp\left\{-\frac{n\rho}{2}\left(\frac{\rho}{V(f)} \wedge 1\right)\right\}. \tag{29}$$

*Let $\rho' \in \mathbb{R}_+$ and $g \in L_+^1(\mu)$ be such that $\mu g = 1$, $g \leq e^{\rho'}$, $f$ and $V(g)$ is finite. If, in addition, $(\rho + \rho') < H(f)$, then the following bound holds true:*

$$\mathrm{E}_g^n(1 - \phi_{n,f,\rho,c}) \leq \exp\left\{-\frac{n[H(f) - (\rho + \rho')]}{2}\left(\frac{H(f) - (\rho + \rho')}{V(g)} \wedge 1\right)\right\}. \tag{30}$$

PROOF. $H(f)$ is also finite. Let us denote $\log f$ by $\ell_f$, $\log g$ by $\ell_g$ and set $s \in (0, 1]$. Then

$$\begin{aligned} c\,\mathrm{E}^{\star n} \phi_{n,f,\rho,c} &= cP^{\star n}\{\ell_{n,f} - \ell_n^\star \geq n\rho - nH(f) + \log c\} \\ &\leq e^{-ns(\rho - H(f))}(P^\star e^{s(\ell_f - \ell^\star)})^n. \end{aligned}$$

A Taylor expansion of the function $t \mapsto P^\star e^{t(\ell_f - \ell^\star)}$ implies that

$$\begin{aligned} P^\star e^{s(\ell_f - \ell^\star)} &= 1 - sH(f) + s^2 \int_0^1 (1 - t) \int (f^\star)^{1 - st} f^{st} (\ell^\star - \ell_f)^2 \, d\mu \, dt \\ &\leq 1 - sH(f) + s^2 V(f^\star, f)^{1 - st} V(f, f^\star)^{st}/2 \\ &\leq 1 - sH(f) + s^2 V(f)/2, \end{aligned}$$

by a Hölder inequality with parameters $1/st$ and $1/(1-st)$. Moreover, since $\log t \le t-1$ $(t>0)$, we have

$$c\,\mathrm{E}^{\star n}\phi_{n,f,\rho,c} \le \exp[-ns\rho + ns^2V(f)/2].$$

The choice $s = 1 \wedge \frac{\rho}{V(f)}$ yields (29). Similarly, for all $s \in (0,1]$,

$$\begin{aligned}\mathrm{E}_g^n(1-\phi_{n,f,\rho,c}) &\le P_g^n\{\ell_n^\star - \ell_{n,f} > n[H(f)-\rho]\}\\ &\le P_g^n\{\ell_n^\star - \ell_{n,g} > n[H(f)-(\rho+\rho')]\}\\ &\le e^{-ns[H(f)-(\rho+\rho')]}(P_g e^{s(\ell^\star-\ell_g)})^n.\end{aligned}$$

The same arguments as before lead to $P_g e^{s(\ell^\star-\ell_g)} \le 1 + s^2V(g)/2$ and

$$\mathrm{E}_g^n(1-\phi_{n,f,\rho,c}) \le \exp\{-ns[H(f)-(\rho+\rho')] + ns^2V(g)/2\}.$$

The choice $s = 1 \wedge \frac{H(f)-(\rho+\rho')}{V(g)}$ yields (30). $\square$

## APPENDIX C: PROOF OF PROPOSITION 1

Let $\{\delta_n\}$ be a decreasing sequence of positive numbers which tend to 0. Let us denote by $\|\cdot\|$ the $L^1(\mu)$ norm. Because $\sqrt{H(\theta)} \ge \|f^\star - f_\theta\|/2$, **M1** ensures that Proposition 1 holds if

$$\pi_{k^\star+1}\{\theta \in \Theta_{k^\star+1} : \|f^\star - f_\theta\| \le \sqrt{\delta_n}\} \le C_2\sqrt{\delta_n}^{D(k^\star)+1} \tag{31}$$

for some $C_2 > 0$ which does not depend on $\{\delta_n\}$. We use a new parameterization for translating $\|f^\star - f_\theta\| \le \sqrt{\delta_n}$ in terms of parameters $\boldsymbol{p}$ and $\boldsymbol{\gamma}$. It is a variant of the locally conic parameterization [6], using the $L^1$ norm instead of the $L^2$ norm. In the sequel, $c, C$ will be generic positive constants.

*$L^1$ locally conic parameterization.* For each $\theta = (\boldsymbol{p},\boldsymbol{\gamma}) \in \mathrm{int}(\Theta_{k^\star+1})$, we define iteratively the permutation $\sigma_\theta$ upon $\{1,\dots,k^\star+1\}$ as follows:

- $(j_1,\sigma_\theta(j_1)) = \min_{(j,j')} \arg\min\{|\gamma_j^\star - \gamma_{j'}|_1 : j \le k^\star, j' \le k^\star+1\}$, where the first minimum is for the lexicographic ranking;
- if $(j_1,\sigma_\theta(j_1)),\dots,(j_{l-1},\sigma_\theta(j_{l-1}))$ with $l < k^\star$ have been defined, then $(j_l, \sigma_\theta(j_l)) = \min_{(j,j')} \arg\min\{|\gamma_j^\star - \gamma_{j'}|_1\}$, where in the arg min, index $j \le k^\star$ does not belong to $\{j_1,\dots,j_{l-1}\}$ and index $j' \le k^\star+1$ does not belong to $\{\sigma_\theta(j_1),\dots,\sigma_\theta(j_{l-1})\}$;
- once $(j_1,\sigma_\theta(j_1)),\dots,(j_{k^\star},\sigma_\theta(j_{k^\star}))$ are defined, the value of $\sigma_\theta(k^\star+1)$ is uniquely determined.

We can assume without loss of generality that $\sigma_\theta = \mathrm{id}$, the identity permutation over $\{1,\dots,k^\star+1\}$. Indeed, for every $\theta = (\boldsymbol{p},\boldsymbol{\gamma}) \in \Theta_{k^\star+1}$ and each permutation $\varsigma$ onto $\{1,\dots,k^\star+1\}$, let $\theta^\varsigma = (\boldsymbol{p}^\varsigma,\boldsymbol{\gamma}^\varsigma) \in \Theta_{k^\star+1}$ be the parameter with coordinates $p_j^\varsigma = p_{\varsigma(j)}$, $\gamma_j^\varsigma = \gamma_{\varsigma(j)}$ (all $j \le k^\star + 1$) and set $\pi_{k^\star+1}^\varsigma(\theta) = \pi_{k^\star+1}(\theta^\varsigma)$.

Since for all $\theta$ and $\varsigma$, $\|f^\star - f_\theta\| = \|f^\star - f_{\theta^\varsigma}\|$,

$$\begin{aligned}\pi_{k^\star+1}\{\theta \in \Theta_{k^\star+1} : \|f^\star - f_\theta\| \le \sqrt{\delta_n}\} \\ = \sum_\varsigma \pi_{k^\star+1}^\varsigma\{\theta \in \Theta_{k^\star+1} : \sigma_\theta = \mathrm{id}, \|f^\star - f_\theta\| \le \sqrt{\delta_n}\},\end{aligned} \tag{32}$$

where the sum above is on all possible permutations.

We show below that the term in the sum above associated with $\varsigma = \mathrm{id}$ is bounded by a constant times $\sqrt{\delta_n}^{D(k^\star)+1}$. The proof involves only properties that all $\pi_{k^\star+1}^\varsigma$ share. Studying the latter term is therefore sufficient to conclude that Proposition 1 holds.

Set $\Theta^\star = \{\theta \in \Theta_{k^\star+1} : \sigma_\theta = \mathrm{id}\}$. For all $\theta \in \Theta^\star$, let $\gamma_\theta = \gamma_{k^\star+1}$, $p_\theta = p_{k^\star+1}$ and $R_\theta = (\rho_1, \ldots, \rho_{k^\star-1}, r_1, \ldots, r_{k^\star})$, where

$$\rho_j = \frac{p_j - p_j^\star}{p_\theta} \quad \text{and} \quad r_j = \frac{\gamma_j - \gamma_j^\star}{p_\theta} \qquad (j \le k^\star).$$

Note that $\sum_{j \le k^\star} \rho_j = -1$. Now, define

$$N(\gamma_\theta, R_\theta) = \left\| g_{\gamma_\theta} + \sum_{j=1}^{k^\star} p_j^\star r_j^T \nabla g_{\gamma_j^\star} + \sum_{j=1}^{k^\star} \rho_j g_{\gamma_j^\star} \right\|,$$

then $t_\theta = p_\theta N(\gamma_\theta, R_\theta)$.

Lemma C.1. *For all $\theta \in \Theta^\star$, let $\Psi(\theta) = (t_\theta, \gamma_\theta, R_\theta)$. The function $\Psi$ is a bijection between $\Theta^\star$ and $\Psi(\Theta^\star)$. Furthermore, $T = \sup_{\theta \in \Theta^\star} t_\theta$ is finite, so that the projection of $\Psi(\Theta^\star)$ along its first coordinate is included in $[0, T]$. Finally, for all $\varepsilon > 0$, there exists $\eta > 0$ such that, for every $\theta \in \Theta^\star$, $\|f^\star - f_\theta\| \le \eta$ yields $t_\theta \le \varepsilon$.*

Proof. It is readily seen that $\Psi$ is a bijection. We point out that $N(\gamma, R)$ is necessarily positive for all $(t, \gamma, R) \in \Psi(\Theta^\star)$, by virtue of **M5**. As for the finiteness of $T$, note that, for any $\theta \in \Theta^\star$,

$$\begin{aligned} t_\theta &= \left\| p_\theta g_{\gamma_\theta} + \sum_{j=1}^{k^\star} p_j^\star (\gamma_j - \gamma_j^\star)^T \nabla g_{\gamma_j^\star} + \sum_{j=1}^{k^\star} (p_j - p_j^\star) g_{\gamma_j^\star} \right\| \\ &\le 2 + \sum_{j=1}^{k^\star} p_j^\star \|(\gamma_j - \gamma_j^\star)^T \nabla g_{\gamma_j^\star}\|. \end{aligned} \tag{33}$$

The right-hand side term above is finite because $\Gamma$ is bounded and $\|(\nabla g_{\gamma_j^\star})_l\|$ $(j \le k^\star, l \le d)$ are finite thanks to **M4**. Hence, $T$ is finite.

The last part of the lemma is a straightforward consequence of the compactness of $\Gamma$ and continuity of $\theta \mapsto f_\theta(z)$. $\square$

*Proof of* (31). For any $\tau > 0$, define the sets

$$B_1^\tau = \left\{\theta \in \Theta^\star : \min_{j \le k^\star} |\gamma_\theta - \gamma_j^\star|_1 > \tau, \|f^\star - f_\theta\| \le \sqrt{\delta_n}\right\}$$

and

$$B_2^\tau = \left\{\theta \in \Theta^\star : \min_{j \le k^\star} |\gamma_\theta - \gamma_j^\star|_1 \le \tau, \|f^\star - f_\theta\| \le \sqrt{\delta_n}\right\}.$$

Inequality (31) is a consequence of the following:

LEMMA C.2. *Given $\tau > 0$, there exists $C > 0$ such that, for all $n \ge 1$,*

$$\pi_{k^\star+1}\{B_1^\tau\} \le C\sqrt{\delta_n}^{k^\star(d+1)}. \tag{34}$$

LEMMA C.3. *There exist $\tau, C > 0$ such that, for all $n \ge 1$,*

$$\pi_{k^\star+1}\{B_2^\tau\} \le C\sqrt{\delta_n}^{k^\star(d+1)}. \tag{35}$$

Because $\Gamma$ is compact, continuity arguments on the norm in finite dimensional spaces yield the following useful property: Under **M5**, if $g_1, \dots, g_k \in L^1(\mu)$ are $k$ functions such that, for every $\gamma \in \Gamma$, $g_\gamma, g_1, \dots, g_k$ are linearly independent, then there exists $C > 0$ such that, for all $a = (a_0, \dots, a_k) \in \mathbb{R}^{k+1}$ and $\gamma \in \Gamma$,

$$\left\|a_0 g_\gamma + \sum_{j=1}^k a_j g_j\right\| \ge C \sum_{j=0}^k |a_j|. \tag{36}$$

PROOF OF LEMMA C.2. Let $\tau > 0$, let $(t, \gamma, R) \in \Psi(\Theta^\star)$ and $\theta = (\boldsymbol{p}, \boldsymbol{\gamma}) = \Psi^{-1}(t, \gamma, R)$ satisfy $|\gamma_\theta - \gamma_j^\star|_1 > \tau$ for all $j \le k^\star$ and $\|f^\star - f_\theta\| \le \sqrt{\delta_n}$. Given any $z \in \mathcal{Z}$, a Taylor–Lagrange expansion (in $t$) of $[f^\star(z) - f_\theta(z)]$ yields the existence of $t^o \in (0, t)$ (depending on $z$) such that

$$\begin{aligned} |f^\star(z) - f_\theta(z)| \ge \frac{t}{N}&\left|g_\gamma(z) + \sum_{j=1}^{k^\star} p_j^\star r_j^T \nabla g_{\gamma_j^\star}(z) + \sum_{j=1}^{k^\star} \rho_j g_{\gamma_j^\star}(z)\right| \\ &- \frac{t^2}{N^2}\left|\sum_{j=1}^{k^\star} \rho_j r_j^T \nabla g_{\gamma_j^o}(z) + \frac{1}{2}\sum_{j=1}^{k^\star} p_j^o r_j^T D^2 g_{\gamma_j^o}(z) r_j\right|, \end{aligned}$$

where $\gamma_j^o = \gamma_j^\star + t^o r_j / N$ and $p_j^o = p_j^\star + t^o \rho_j / N$ (all $j \le k^\star$). Therefore, by virtue of **M4**, there exists $C > 0$ such that

$$\|f^\star - f_\theta\| \ge t\left(1 - C\frac{t}{N^2}\left[\sum_{j=1}^{k^\star} (|\rho_j||r_j|_1 + |r_j|_2^2)\right]\right). \tag{37}$$

Furthermore, **M5** and (36) imply that, for some $C>0$ (depending on $\tau$),

$$N \geq C\left(1+\sum_{j=1}^{k^\star}(|\rho_j|+p_j^\star|r_j|_1)\right), \tag{38}$$

so the following lower bound on $\|f^\star - f_\theta\|$ is deduced from (37):

$$\|f^\star - f_\theta\| \geq t\left(1 - C\frac{\sum_{j=1}^{k^\star}(|p_j-p_j^\star||\gamma_j-\gamma_j^\star|_1 + |\gamma_j-\gamma_j^\star|_2^2)}{\sum_{j=1}^{k^\star}(|p_j-p_j^\star|+p_j^\star|\gamma_j-\gamma_j^\star|_1)}\right). \tag{39}$$

By mimicking the last part of the proof of Lemma C.1, we obtain that the right-hand term in (39) is larger than $t/2$ for $n$ large enough (independently of $\theta$). Because $t=p_\theta N$ and (38) holds, there exists $c>0$ such that

$$\pi_{k^\star+1}\{B_1^\tau\} \leq \pi_{k^\star+1}\left\{\theta\in\Theta^\star : \sum_{j=1}^{k^\star}(|p_j-p_j^\star|+p_j^\star|\gamma_j-\gamma_j^\star|_1) \leq c\sqrt{\delta_n}\right\},$$

leading to (34) by virtue of **M1**. □

Proof of Lemma C.3. Let $\tau>0$ and $\theta=(\boldsymbol{p},\boldsymbol{\gamma})\in\Theta^\star$ satisfying $\|f^\star - f_\theta\|\leq\sqrt{\delta_n}$. Assume that $|\gamma_\theta-\gamma_j^\star|_1\leq\tau$ for some $j\leq k^\star$, say, $j=1$. By construction of $\Theta^\star$, $|\gamma_1-\gamma_1^\star|_1\leq|\gamma_\theta-\gamma_1^\star|_1\leq\tau$, and $\tau$ can be chosen small enough so that $\gamma_\theta$ must be different from $\gamma_j^\star$ for every $j=2,\dots,k^\star$. We consider without loss of generality that $\gamma_\theta\notin\{\gamma_j^\star : j\leq k^\star\}$.

Lemma C.1 implies that $|\gamma_j-\gamma_j^\star|_1$ and $|p_j-p_j^\star|$ go to 0 as $n\uparrow\infty$ for every $j=2,\dots,k^\star$. This yields that $|p_1+p_\theta-p_1^\star|$ goes to 0 as $n\uparrow\infty$. Therefore, by virtue of **M5** and (36), there exist $c,C>0$ such that, for $n$ large enough,

$$\begin{aligned}\|f^\star - f_\theta\| \geq C\Bigg(&\sum_{j=2}^{k^\star}|p_j-p_j^\star| + \sum_{j=2}^{k^\star}p_j^\star|\gamma_j-\gamma_j^\star|_1\\ &+\Bigg|(p_1+p_\theta-p_1^\star)\\ &\quad+\sum_{(r,s)\notin\mathcal{A}}\lambda_{rs}^0[p_\theta(\gamma_\theta-\gamma_1^\star)_r(\gamma_\theta-\gamma_1^\star)_s\\ &\qquad\qquad+p_1(\gamma_1-\gamma_1^\star)_r(\gamma_1-\gamma_1^\star)_s]\Bigg|\\ &+\sum_{(r,s)\in\mathcal{A}}[p_\theta|(\gamma_\theta-\gamma_1^\star)_r(\gamma_\theta-\gamma_1^\star)_s|\\ &\qquad+p_1|(\gamma_1-\gamma_1^\star)_r(\gamma_1-\gamma_1^\star)_s|]\end{aligned} \tag{40}$$

$$+\sum_{l=1}^{d}\Bigg|p_1(\gamma_1-\gamma_1^\star)_l+p_\theta(\gamma_\theta-\gamma_1^\star)_l$$
$$+\sum_{(r,s)\notin\mathcal{A}}\lambda_{rs}^l[p_\theta(\gamma_\theta-\gamma_1^\star)_r(\gamma_\theta-\gamma_1^\star)_s$$
$$+p_1(\gamma_1-\gamma_1^\star)_r(\gamma_1-\gamma_1^\star)_s]\Bigg|\Bigg)$$
$$-c\Bigg(p_\theta|\gamma_\theta-\gamma_1^\star|_1^3+p_1|\gamma_1-\gamma_1^\star|_1^3+\sum_{j=2}^{k^\star}|\gamma_j-\gamma_j^\star|_2^2\Bigg)$$
$$=CA_1-cA_2.$$

Since $|\gamma_j-\gamma_j^\star|_1$ goes to 0 for $j=2,\dots,k^\star$, $\sum_{j=2}^{k^\star}|\gamma_j-\gamma_j^\star|_2^2$ can be neglected compared to $\sum_{j=2}^{k^\star}p_j^\star|\gamma_j-\gamma_j^\star|_1$ when $n$ is large enough. If $CA_1\le 2cA_2$, then $\sum_{j=2}^{k^\star}|p_j-p_j^\star|\le 2cA_2$, so that $|p_1+p_\theta-p_1^\star|\le 2cA_2$, which yields in turn

$$\Bigg|\sum_{(r,s)\notin\mathcal{A}}\lambda_{rs}^0[p_\theta(\gamma_\theta-\gamma_1^\star)_r(\gamma_\theta-\gamma_1^\star)_s+p_1(\gamma_1-\gamma_1^\star)_r(\gamma_1-\gamma_1^\star)_s]\Bigg|$$
$$+\sum_{(r,s)\in\mathcal{A}}[p_\theta|(\gamma_\theta-\gamma_1^\star)_r(\gamma_\theta-\gamma_1^\star)_s|+p_1|(\gamma_1-\gamma_1^\star)_r(\gamma_1-\gamma_1^\star)_s|]\le 4cA_2.$$

Consequently, **M5** guarantees the existence of $C'>0$ such that

$$p_\theta|\gamma_\theta-\gamma_1^\star|_2^2+p_1|\gamma_1-\gamma_1^\star|_2^2$$
$$\le C'(p_\theta|\gamma_\theta-\gamma_1^\star|_1^3+p_1|\gamma_1-\gamma_1^\star|_1^3),$$

which is impossible when $\tau$ is chosen small enough. Therefore, $CA_1>2cA_2$ and (40) together with **M5** give

$$\|f^\star-f_\theta\|\ge C\Bigg(\sum_{j=2}^{k^\star}|p_j-p_j^\star|+\sum_{j=2}^{k^\star}p_j^\star|\gamma_j-\gamma_j^\star|_1+|p_1+p_\theta-p_1^\star|$$
$$+p_\theta|\gamma_\theta-\gamma_1^\star|_2^2+p_1|\gamma_1-\gamma_1^\star|_2^2$$
$$+\sum_{l=1}^{d}\Bigg|p_1(\gamma_1-\gamma_1^\star)_l+p_\theta(\gamma_\theta-\gamma_1^\star)_l$$
$$+\sum_{(r,s)\notin\mathcal{A}}\lambda_{rs}^l[p_\theta(\gamma_\theta-\gamma_1^\star)_r(\gamma_\theta-\gamma_1^\star)_s$$
$$+p_1(\gamma_1-\gamma_1^\star)_r(\gamma_1-\gamma_1^\star)_s]\Bigg|\Bigg),$$

for some $C > 0$. Finally,

$$
\begin{aligned}
&|p_1 + p_\theta - p_1^\star| + \sum_{j=2}^{k^\star} |p_j - p_j^\star| + p_1|\gamma_1 - \gamma_1^\star|_2^2 + p_\theta|\gamma_\theta - \gamma_1^\star|_2^2 \\
&\qquad + |p_1(\gamma_1 - \gamma_1^\star) + p_\theta(\gamma_\theta - \gamma_1^\star)|_1 + \sum_{j=2}^{k^\star} p_j^\star|\gamma_j - \gamma_j^\star|_1 \le C\sqrt{\delta_n}.
\end{aligned} \tag{41}
$$

Therefore, for $\tau$ small enough and $n$ large enough,

$$\pi_{k^\star+1}\{B_2^\tau\} \le \pi_{k^\star+1}\{\theta \in \Theta^\star : (41) \text{ holds}\}.$$

The conditions on $p_j$ and $\gamma_j$ $(j = 2, \ldots, k^\star)$ and a symmetry argument imply that the right-hand side term above is bounded by a constant times $\sqrt{\delta_n}^{[(d+1)(k^\star-1)]}$ times $w_n$, where

$$
\begin{aligned}
w_n = \int \mathbb{1}\{p_\theta \ge p_1\}\mathbb{1}\{&|p_1 + p_\theta - p_1^\star| + p_1|\gamma_1 - \gamma_1^\star|_2^2 \\
&+ p_\theta|\gamma_\theta - \gamma_1^\star|_2^2 + |p_1(\gamma_1 - \gamma_1^\star) \\
&+ p_\theta(\gamma_\theta - \gamma_1^\star)|_1 \le C\sqrt{\delta_n}\}\, d\pi^{\boldsymbol{\gamma}}_{k^\star+1}(\boldsymbol{\gamma})\, d\pi^{\boldsymbol{p}}_{k^\star+1}(\boldsymbol{p}).
\end{aligned}
$$

Note that the conditions in the integrand imply that $|\gamma_\theta - \gamma_1|_2^2 \le 4C\sqrt{\delta_n}/p_1$ and $p_\theta \ge p_1^\star/4$ as soon as $C\sqrt{\delta_n} \le p_1^\star/2$. Simple calculus (based on **M1**) yields the result.

## APPENDIX D: PROOF OF PROPOSITION 2

It is readily seen that $\mathbf{O2}(k^\star + 1)$ holds for the chosen approximating set. Let us focus now on the entropy condition (6).

*Constructing $\delta$-brackets.* Let $\delta_1$ satisfy (5). A convenient value will be chosen later on. Set $j' \le \lfloor \delta_0/\delta_n \rfloor$, $\varepsilon = j'\delta_n/4$ and $\tau \ge 1$.

Let $\theta = (\boldsymbol{p}, \boldsymbol{\gamma}) \in \Theta_{k^\star+1}$ be arbitrarily chosen. Let $\eta \in (0, \eta_1)$ be small enough so that, for every $j \le k^\star + 1$, $u_j = \overline{g}_{\gamma_j,\eta}$ and $v_j = \underline{g}_{\gamma_j,\eta}$ (as defined in **M2**) satisfy, for all $\gamma \in \Gamma$, $P_{u_j - v_j}(1 + \log^2 g_\gamma) \le \varepsilon/\tau$, $P_{g_\gamma}(\log u_j - \log v_j)^2 \le (\varepsilon/\tau)^2$ and $P_{u_j - v_j} \log^2 u_j \le (\varepsilon/\tau)\log^2(\varepsilon/\tau)$.

If we define $v_\theta = (1 - \varepsilon/\tau)(\sum_{j=1}^{k^\star+1} p_j v_j)$ and $u_\theta = (1 + \varepsilon/\tau)(\sum_{j=1}^{k^\star+1} p_j u_j)$, then there exists $\tau \ge 1$ (which depends only on $k^\star$ and the constant $M$ of **M2**) such that the bracket $[v_\theta, u_\theta]$ is an $\varepsilon$-bracket. The repeated use of $(\sum_j p_j u_j / \sum_j p_j v_j) \le \max_j u_j/v_j$ is the core of the proof we omit.

*Control of the entropy.* The rule $x_1(1-\varepsilon/\tau) = e^{-n}$ and $x_{j+1}(1-\varepsilon/\tau) = x_j(1+\varepsilon/\tau)$ is used for defining a net for the interval $(e^{-n}, 1)$. Such a net has at most $[1 + n/\log(1+\varepsilon/\tau)/(1-\varepsilon/\tau) \le 1 + 2n\tau/\varepsilon]$ support points. Using repeatedly this construction on each dimension of the $(k^\star+1)$-dimensional simplex yields a net for $\{\boldsymbol{p} \in \mathbb{R}_+^{k^\star} : \min_{j \le k^\star} p_j \ge e^{-n}, 1 - \sum_{j \le k^\star} p_j \ge e^{-n}\}$ with at most $O((n/\varepsilon)^{(k^\star+1)})$ support points. We can choose a net for $\Gamma^{k^\star+1}$ with at most $O(\varepsilon^{-d(k^\star+1)})$ support points such that each $\boldsymbol{\gamma} \in \Gamma^{k^\star+1}$ is within $|\cdot|_1$-distance $\varepsilon$ of some element of the net.

Consequently, the minimum number of $\varepsilon$-brackets needed to cover $\mathcal{F}_n^{k^\star+1}$ is $O(n^{(k^\star+1)}/\varepsilon^{(d+1)(k^\star+1)})$, so there exist constants $a, b, c > 0$ for which

$$\mathcal{E}\left(\mathcal{F}_n^{k^\star+1}, \frac{j'\delta_n}{4}\right) \le a \log n - b \log(j'\delta_n) + c. \tag{42}$$

Now, let us note that $\frac{nj'\delta_n}{\log^2(j'\delta_n)} \ge \frac{n\delta_n}{(\log\delta_n)\log(j'\delta_n)} \ge \frac{n\delta_n}{\log^2\delta_n}$ and consider each term of the right-hand side of (42) in turn. It is readily seen that $a\log n \le n\delta_n/\log^2\delta_n$ is equivalent to

$$\delta_1 \ge [(\log^3 n) n^{(\delta_1/a)^{1/2}-1}]^{-1}. \tag{43}$$

Now, $-b\log(j'\delta_n) \le \frac{n\delta_n}{(\log\delta_n)\log(j'\delta_n)}$ if and only if $-b\log\delta_n \le \delta_1 \log^3 n$. Since $\log^2\delta_n \le 4\log^2 n$, both are valid as soon as

$$\delta_1 \ge 2b/\log^2 n. \tag{44}$$

Finally, using again $\log^2\delta_n \le 4\log^2 n$ yields that $c \le n\delta_n/\log^2\delta_n$ when

$$\delta_1 \ge 4c/\log n. \tag{45}$$

When $\delta_1 \ge a$, the largest values of the right-hand sides of (43), (44) and (45) are achieved at $n_0$. So, $\delta_1$ can be chosen large enough (independently of $j'$ and $n$) so that (5), (43), (44) and (45) hold for all $n \ge n_0$ and $j' \le \lfloor \delta_0/\delta_n \rfloor$. This completes the proof of Proposition 2, because $\mathcal{E}(\mathcal{F}_n^{k^\star+1}, j'\delta_n/4)$ is larger than the left-hand side of (6) (with $j'$ substituted to $j$). □

**Acknowledgments.** We thank the referees and Associate Editor for their helpful suggestions.

MAP5 CNRS UMR 8145
Université Paris Descartes
45 rue des Saints-Pères
75270 Paris cedex 06
France
E-mail: chambaz@univ-paris5.fr

Cérémade, UMR CNRS 7534
Université Paris Dauphine and CREST
Place de Lattre de Tassigny
75775 Paris cedex 16
France
E-mail: rousseau@ceremade.dauphine.fr